

\baselineskip=14pt
\parskip=10pt

\magnification=\magstephalf

\def\1{{\overline{1}}}
\def\2{{\overline{2}}}
\parindent=0pt
\overfullrule=0in

\def\frac#1#2{{#1 \over #2}}
\centerline
{\bf Tweaking the Beukers Integrals In Search of More Miraculous Irrationality Proofs \`A La Ap\'ery }
\bigskip
\centerline
{\it Robert DOUGHERTY-BLISS, Christoph KOUTSCHAN, and Doron ZEILBERGER}

\bigskip

\qquad \qquad {\it In honor of our irrational guru Wadim Zudilin, on his $\lfloor 50 \, \zeta(5) \rfloor$-th birthday}

{\bf [Actual] Historical Introduction: How Beukers' Proofs Were ACTUALLY found}

{\bf  Hilbert's 0-th problem} 

Before David Hilbert [H] stated his famous 23 problems, he mentioned two problems
that he probably believed to be yet much harder, and indeed, are still wide open today. One of them was to prove
that there are infinitely many prime numbers of the form $2^n+1$, and the other one was to prove that the
Euler-Mascheroni constant is irrational.

Two paragraphs later he stated his optimistic belief that ``in mathematics there is no {\it ignorabimus}.''

As we all know, he was proven wrong by G\"odel and Turing
{\it in general}, but even
for such  concrete  problems, like the irrationality of a specific, natural, constant, like  the Euler-Mascheroni constant
(that may be defined in terms of the {\bf definite} integral $\quad -\int_0^{\infty} e^{-x} \log x$) , that is most probably
{\bf decidable} in the logical sense, (i.e. there {\it probably} {\bf exists} a (rigorous) proof), 
we {\it lowly humans} did not yet find it, (and may never will!).

While the Euler-Mascheroni constant (and any other, natural,  explicitly-defined, constant that is not obviously rational) is {\bf surely} 
irrational, in the everyday sense of the word {\it sure} (like death and taxes), giving a {\it proof}, in the {\bf mathematical} sense of `proof'
is a different matter. While $e$ was proved irrational a long time ago (trivial exercise), and $\pi$ was proved irrational
by Lambert around 1750, we have no clue how to prove that $e+\pi$ is irrational.
Ditto for  $e \cdot \pi$. {\bf Exercise}: Prove that {\it at least} one of them is irrational.

{\bf Ap\'ery's Miracle}

As Lindemann first proved in 1882,
the number $\pi$ is more than {\it just} irrational, it is {\it transcendental}, hence it follows that $\zeta(n)$ is irrational
for all even arguments, since Euler proved that $\zeta(2n)$ is a multiple of $\pi^{2n}$ by a rational number. But
proving that $\zeta(3)$, $\zeta(5)$, $\dots$ are irrational remained wide open.

Since such problems are so hard, it was breaking news, back in 1978, when 64-years-old Roger Ap\'ery announced and sketched
a proof that $\zeta(3):=\sum_{n=1}^{\infty} \frac{1}{n^3}$ is irrational. This was beautifully narrated in a classic 
expository paper by Alf van der Poorten [vdP], aided by details filled-in by Henri Cohen and Don Zagier.
While beautiful in our eyes, most people found the proof ad-hoc and too complicated, and they did not like
the heavy reliance on {\it recurrence relations}.

To those people, who found Ap\'ery's original proof too {\it magical}, ad-hoc, and {\it computational}, another proof,
by a 24-year-old PhD student by the name of Frits Beukers [B] was a breath of fresh air. It was a marvelous {\it gem}
in human-generated mathematics, and could be easily followed by a first-year student, using {\it partial fractions} and
very easy estimates of a certain triple integral, namely

$$
\int_0^1 \, \int_0^1 \, \int_0^1 \, \frac{(x(1-x)y(1-y)z(1-z))^n}{(1-z+xyz)^{n+1}} \,dx\,dy\,dz\quad .
$$

The general approach of Ap\'ery of finding  concrete sequences of integers $a_n,b_n$ such that 
$$
|\zeta(3) -\frac{a_n}{b_n}| \, < \, \frac{CONST}{b_n^{1+\delta}} \quad,
$$
(see below) for a {\bf positive} $\delta$ was still followed, but the details were much more palatable and {\it elegant} to the average 
{\it mathematician in the street}.

As a warmup, Beukers, like Ap\'ery before him, gave a new proof of the already proved fact that $\zeta(2)=\frac{\pi^2}{6}$ is irrational, using
the double integral
$$
 \int_0^1 \, \int_0^1 \, \frac{(x(1-x)y(1-y))^n}{(1-xy)^{n+1}} \, dx \, dy \quad .
$$

Ironically, we will follow Beukers' lead, but heavily using {\bf recurrence relations}, that will be the {\bf engine}
of our approach. Thus we will abandon the original {\it raison d'\^etre} of Beukers' proof of 
{\bf getting rid of recurrences}, and bring them back with a {\bf vengeance}.

{\bf [Alternative World] Historical Introduction: How Beukers's Proofs Could (and Should!) have been Discovered}

Once upon a time, there was a precocious teenager, 
who was also a {\it computer whiz}, let's call him/her/it/they {\it Alex}. Alex just got a new laptop
that had {\tt Maple},  as a birthday present.

Alex typed, {\it for no particular reason},

{\tt int(int(1/(1-x*y),x=0..1),y=0..1);}

and immediately got the answer: $\frac{\pi^2}{6}$. Then Alex was wondering about the sequence
$$
I(n):=  \int_0^1 \, \int_0^1 \, \frac{(x(1-x)y(1-y))^n}{(1-xy)^{n+1}} \, dx \, dy \quad .
$$
(why not, isn't it a natural thing to try out for a curious teenager?), and typed

{\tt I1:=n->int(int(1/(1-x*y)*(x*(1-x)*y*(1-y)/(1-x*y))**n,x=0..1),y=0..1);}

({\tt I} is reserved in Maple for $\sqrt{-1}$, so Alex needed to use {\tt I1}),

and looked at the first ten values by typing:

{\tt L:=[seq(I1(i),i=1..10)];} \quad ,

getting after a few seconds

$$
[5- \frac{\pi^{2}}{2} , -{\frac{125}{4}}+{\frac {19\,{\pi}^{2}}{6}},{\frac{8705}{36}}-{\frac {49\,{\pi}^{2}}{2}},-{\frac{32925}{16}}+{\frac {417\,{\pi}^{2}}{2}},
$$
$$
{\frac{13327519}{720}}-{\frac {3751\,{\pi}^{2}}{2}},-{\frac{124308457}{720}}+{\frac {104959\,{\pi}^{2}}{6}},
$$
$$
{\frac{19427741063}{11760}}-{\frac {334769\,{\pi}
^{2}}{2}},-{\frac{2273486234953}{141120}}+{\frac {9793891\,{\pi}^{2}}{6}},
$$
$$
{\frac{202482451324891}{1270080}}-{\frac {32306251\,{\pi}^{2}}{2}},-{\frac{
2758128511985}{1728}}+{\frac {323445423\,{\pi}^{2}}{2}}] \quad .
$$

Alex immediately noticed that, at least for $n\leq 10$,  
$$
I(n)= a_n - b_n \frac{\pi^2}{6} \quad ,
$$ 
for some {\bf integers} $b_n$ and some
{\bf rational numbers} $a_n$.  By taking {\tt evalf(L)}, Alex also noticed that $I(n)$ get smaller and smaller. Knowing that Maple
could not be trusted with floating point calculations (unless you change the  value of {\tt Digits} from its default, to something
higher, say, in this case {\tt Digits:=30}), that they get smaller and smaller. Typing `{\tt evalf(L,30)};',  Alex got:
$$
[ 0.06519779945532069058275450006, 0.0037472701163022929758881663, 
$$
$$
0.000247728866269394110526059, 0.00001762713127202699137347, 
$$
$$
0.0000013124634659314676853, 0.000000100776323486001254, 
$$
$$
0.00000000791212964371946, 0.0000000006317437711206,
$$
$$
{ 5.1111100706\times 10^{-11}},{ 4.17922459\times 10^{-12}}] \quad.
$$

Alex realized that $I(n)$ seems to go to $zero$ fairly fast,
and since $I(10)/I(9)$ and $I(9)/I(8)$ were pretty close, Alex conjectured that the limit
of $I(n)/I(n-1)$ tends to a certain constant. But ten data points do not suffice!

When Alex tried to find the first $2000$ terms, Maple got slower and slower. Then Alex asked Alexa, the famous robot,

{\it Alexa: how do I compute many terms of the sequence $I(n)$ given by that double-integral?}

and Alexa replied:

{\it Go to Doron Zeilberger's web-site and download the amazing program}

{\tt https://sites.math.rutgers.edu/\~{}zeilberg/tokhniot/MultiAlmkvistZeilberger.txt} \quad,

that accompanied the article [ApaZ]. Typing

{\tt MAZ(1,1/(1-x*y),x*(1-x)*y*(1-y)/(1-x*y),[x,y],n,N, $\{\}$)[1];}

immediately gave a recurrence satisfied by $I(n)$

$$
I(n)=
-{\frac { \left( 11\,{n}^{2}-11\,n+3 \right) }{{n}^{2}}} \cdot {\it I} \left( n-1 \right)
+{\frac { \left( n-1 \right) ^{2} }{{n}^{2}}}\cdot{\it I} \left( n-2 \right) \quad .
$$

Using this recurrence, Alex easily computed the first $2000$ terms, using the following Maple one-liner (calling the sequence
defined by the recurrence {\tt I2(n)}):

{\tt I2:=proc(n) option remember: if n=0 then Pi**2/6 elif n=1 then 5-Pi**2/2 else -(11*n**2-11*n+3)/n**2*I2(n-1)+(n-1)**2/n**2*I2(n-2):fi: end:}

and found out that indeed $I(n)/I(n-1)$ tends to a limit, about $0.09016994$.  Writing
$$
I(n)= a_n - b_n \frac{\pi^2}{6} \quad
$$
and realizing that $I(n)$ is small, Alex found terrific {\bf rational approximations} to $\frac{\pi^2}{6}$, $a_n/b_n$, that after
{\bf clearing denominators} can be written as $a'_n/b'_n$ where {\it now} {\bf both}  numerator $a'_n$ and denominator $b'_n$ are {\bf integers}.
$$
\frac{\pi^2}{6} \approx \frac{a'_n}{b'_n} \quad .
$$

Alex also noticed that for all $n$  up to $2000$, for some constant $C$,
$$
|\frac{\pi^2}{6} - \frac{a'_n}{b'_n}| \leq \frac{C}{(b'_n)^{1+\delta}}\quad ,
$$
where $\delta$ is roughly  $0.09215925467$. Then Alex concluded that this proves that $\frac{\pi^2}{6}$ is irrational, since if
it were rational the left side would have been $\geq \frac{C_1}{b'_n}$, for some constant $C_1$. Of course, some details would still need to be filled-in,
but that was not too hard.

{\bf The General Strategy}

Let's follow Alex's lead. (Of course our fictional Alex owes a lot to the real Beukers and also to Alladi and Robinson [AR]).

Start with a constant, let's call it $C$, given by an explicit integral
$$
\int_{0}^{1} K(x) \, dx \quad,
$$
for some {\bf integrand} $K(x)$, or, more generally, a $d$-dimensional integral
$$
\int_{0}^{1} \dots \int_0^1 K(x_1, \dots, x_k) \, dx_1 \dots dx_k \quad .
$$

Our goal in life is to prove that $C$ is irrational. Of course $C$ may turn out to be rational (that happens!),
or more likely, an algebraic number, or expressible in terms
of a logarithm of an algebraic number, for which, there already exist irrationality proofs (albeit not always {\it effective} ones).
But who knows? Maybe this constant has never been proved irrational, and if it will happen to be famous (e.g. Catalan's constant, or
$\zeta(5)$,  or the Euler-Mascheroni constant mentioned above), we will be famous too. But even if it is a nameless constant, it is still
extremely interesting,  if it is the {\it first} irrationality proof, since these proofs are so hard, witness that, in spite
of great efforts by experts like Wadim Zudilin, the proofs of these are still wide open.

In this article we will present numerous candidates. Our proofs of irrationality  are modulo a `divisibility lemma' (see below),
that we are {\it sure} that someone like  Wadim Zudilin, to whom this paper is dedicated, can fill-in. Our only
doubts are whether these constants are not already proved to be irrational because they happen to be algebraic
(probably not, since Maple was unable to {\it identify} them),  or more complicated numbers (like logarithms of algebraic numbers).
Recall that Maple's {\tt identify} can't (yet) identify {\it everything} that God can.

Following Beukers and Alladi-Robinson, we introduce a sequence of integrals, parameterized by a non-negative integer $n$

$$
I(n)=\int_{0}^{1} K(x)\,(x(1-x)K(x))^n \, dx \quad ,
$$
and analogously for multiple integrals, or more generally
$$
I(n)=\int_{0}^{1} K(x)\,(x(1-x)S(x))^n \, dx \quad ,
$$
for another function $S(x)$. Of course $I(0)=C$, our constant that we want to prove irrational.

It so happens that for a wide class of functions $K(x)$, $S(x)$, (for single or multivariable $x$) using 
the {\bf Holonomic ansatz} [Ze1], and implemented  (for the single-variable case) in [AlZ], and for the multi-variable
case in [ApZ], and {\it much more efficiently} in [K], there exists a {\bf linear recurrence equation with polynomial coefficients},
that can be actually computed (always in theory, but also often in practice, unless the dimension is high). In other words
we can find a positive integer $L$, the {\bf order} of the recurrence, and {\bf polynomials} $p_0(n), p_1(n), \dots, p_L(n)$, such that
$$
p_0(n)I(n)+ p_1(n) I(n+1)+ \dots + p_L(n) I(n+L) \, = \, 0 \quad .
$$

If we are lucky (and all the cases in this paper fall into this case) the order $L$ is $2$. 
Furthermore, it would be evident in all the examples in this paper that $p_0(n),p_1(n),p_2(n)$ can be taken to have
{\bf integer coefficients}.

Another `miracle' that happens in all
the examples in this paper is that $I(0)$ and $I(1)$ are rationally-related, i.e. there exist integers $c_0,c_1,c_2$ such that
$$
c_0 I(0)+ c_1 I(1)= c_2 \quad,
$$
that our computers can easily find.

It then follows, by induction, that one can write
$$
I(n)= b_n C  -a_n \quad ,
$$
for some sequences of {\bf rational numbers} $\{a_n\}$ and $\{b_n\}$ that both satisfy the same recurrence as $I(n)$.

Either using trivial bounds on the integral, or using the so-called Poincar\'e lemma (see, e.g. [vdP], [ZeZu1],[ZeZu2]) it turns out that
$$
a_n \,= \,\Omega (\alpha^n) \quad ,  \quad b_n \, = \,\Omega (\alpha^n) \quad , 
$$
for some constant $\alpha>1$, and

$$
|I(n)| = \Omega(\, \frac{1}{\beta^n}\, ) \quad ,
$$
for some constant $\beta>1$.

[Please note that we use $\Omega$ in a looser-than-usual sense, for us $x(n)=\Omega(\alpha^n)$ means that
$\lim_{ n \rightarrow \infty}\, \frac{\log x(n)}{n}=\alpha$.]

In the tweaks of Beukers' integrals for $\zeta(2)$ and $\zeta(3)$  coming up later, $\alpha$ and $\beta$ are equal, but in
the tweaks of the Alladi-Robinson  integrals, $\alpha$ is usually different than $\beta$.

It follows that
$$
|C-\frac{a_n}{b_n}|= \Omega(\frac{1}{(\alpha \beta)^n})  \quad.
$$

Note that $a_n$, and $b_n$ are, usually, {\bf not} integers, but rather {\it rational numbers} 
(In the original Beukers/Ap\'ery cases, the $b_n$  were 
integers, but the $a_n$ were not, in the more general cases in this article, usually neither of them are integers).

It so happens, in all the cases that we discovered, that there exists another sequence of rational numbers $E(n)$ such that
$$
a'_n:=a_n\,E(n) \quad, \quad  b'_n:=b_n\,E(n) \quad ,
$$
are always  {\bf integers}, and, of course $gcd(a'_n \,, \, b'_n)=1$. We call $E(n)$ the {\bf integer-ating factor}.

In some cases we were able to conjecture $E(n)$ exactly,
in terms of products of primes satisfying certain conditions (see below), 
but in other cases we can only conjecture that such an explicitly-describable  sequence exists.

In either  case there exists a real number, that sometimes can be described exactly, and other times only estimated, let's call it $\nu$, such that
$$
\lim_{n \rightarrow \infty}  \frac{\log E(n)}{n} \, = \, \nu \quad ,
$$
or, in our notation, $E(n)=\Omega(\, e^{n \nu} \, )$ .

Since we have
$$
|C-\frac{a'_n}{b'_n}|= \Omega(\frac{1}{(\alpha \beta)^n})  \quad ,
$$
where $b'_n=\Omega(e^{\nu\,n} \alpha^n)$. We need a {\bf positive} $\delta$ such that
$$
(e^{\nu\, n} \alpha^n)^{1+ \delta}=(\alpha \beta)^n \quad.
$$
Taking $\log$ (and dividing by $n$) we have
$$
(\nu + \log \alpha)(1+\delta)= \log \alpha + \log \beta  \quad ,
$$
giving
$$
\delta=\frac{\log \beta - \nu}{\log \alpha + \nu} \quad .
$$
If we are {\it lucky}, and $\log \beta > \nu$, then we have $\delta>0$, and an irrationality proof!, Yea!
We also, at the same time, determined an {\it irrationality measure} (see [vdP])
$$
1+ \frac{1}{\delta} \, = \, \frac{\log \alpha + \log \beta}{\log \beta - \nu}  \quad .
$$
If we are {\it unlucky}, and $\delta<0$, it is still an exponentially fast way to compute our constant $C$ to any desired accuracy.

{\bf Summarizing}: For each specific constant defined by a definite integral, we need to exhibit 

$\bullet$ A second-oder recurrence equation for the numerator and denominator sequence $a_n$ and $b_n$ that feature in
$I(n)=b_n C -a_n$.

$\bullet$ The {\it initial conditions} $a_0,a_1$, $b_0,b_1$ enabling a very fast computation of many terms of $a_n,b_n$.

$\bullet$ The constants $\alpha$ and $\beta$

$\bullet$ Exhibit a conjectured {\bf integer-ating factor} $E(n)$, or else conjecture that one exists, and find,
or estimate (respectively), $\nu := \, \lim_{n \rightarrow \infty} \frac{\log E(n)}{n}$ .

$\bullet$ Verify that $\beta>e^{\nu}$ and get (potentially) famous.

{\bf The Three Classical Cases}

${\bf \log 2}$ ([AR])

$$
C\,= \, \int_0^1 \, \frac{1}{1+x} \, dx \, = \, \log 2 \quad .
$$ 

$$
I(n)= \int_0^1 \frac{(x(1-x))^n}{(1+x)^{n+1}} \, dx \quad .
$$

Recurrence:
$$
\left( n+1 \right) X \left( n \right) + \left( -6\,n-9 \right) X \left( n+1 \right) + \left( n+2 \right) X \left( n+2 \right) \, = 0 \, \quad .
$$

$$
\alpha \,= \, \beta=3+ 2\sqrt{2} \quad .
$$
Initial conditions 
$$
a_0=0 \,, \, a_1=2 \quad ; \quad b_0=1 \, , \, b_1=3 \quad .
$$

{\bf Integer-ating factor} $E(n)=lcm(1 \dots n)$, $\nu=1$.

$$
\delta= \frac{\log \beta -\nu}{\log \alpha + \nu} =
\frac{\log \beta -1}{\log \alpha + 1} =
\frac{\log (3+2\sqrt{2}) -1}{\log (3+ 2\sqrt{2}) + 1} =
0.276082871862633587 \quad.
$$

Implied irrationality measure: $1+1/\delta=4.622100832454231334\dots$.

$\bf{\zeta(2)}$ ([B])

$$
C= \int_0^1\, \int_0^1 \, \frac{1}{1-xy} \, dx \, dy \, = \,\zeta(2) \quad .
$$

$$
I(n)= \int_0^1\, \int_0^1 \, \frac{(x(1-x)y(1-y))^n}{(1-xy)^{n+1}} \,dx \, dy \quad .
$$

Recurrence:
$$
- \left( 1+n \right) ^{2}X \left( n \right) + \left( 11\,{n}^{2}+33\,n+25 \right) X \left( n+1 \right) + \left( 2+n \right) ^{2}X \left( n+2 \right) 
\, = 0 \quad .
$$

$$
\alpha \, =\, \beta= \frac{11}{2} +\frac{5 \sqrt {5}}{2} \quad .
$$
Initial conditions 
$$
a_0=0 \,, \, a_1=-5 \quad ; \quad b_0=1 \, , \, b_1=-3 \quad .
$$

{\bf Integer-ating factor} $E(n)=lcm(1 \dots n)^2$, $\nu=2$.

$$
\delta= \frac{\log \beta -\nu}{\log \alpha + \nu} =
\frac{\log \beta -2}{\log \alpha + 2} =
\frac{\log (11/2+5\sqrt{5}/2) -2}{\log (11/2+5\sqrt{5}/2) + 2} =
 0.09215925473323\dots \quad.
$$

Implied irrationality measure: $1+1/\delta= 11.8507821910523426959528\dots$.

$\bf{\zeta(3)}$ ([B])

$$
C= \int_0^1 \, \int_0^1\, \int_0^1 \, \frac{1}{1-z+xyz} \, dx \, dy \,\,dz = \,\zeta(3) \quad .
$$

$$
I(n)=  \int_0^1 \,\int_0^1\, \int_0^1 \, \frac{(x(1-x) y(1-y) z(1-z))^n}{(1-z+xyz)^{n+1}} \,dx \, dy\,dz \quad .
$$

Recurrence:
$$
\left( 1+n \right) ^{3}X \left( n \right) - \left( 2\,n+3 \right)  \left( 17\,{n}^{2}+51\,n+39 \right) X \left( n+1 \right) + \left( n+2 \right) ^{3}X
 \left( n+2 \right) \, = 0 \quad .
$$

$$
\alpha \, =\, \beta=  17+12\,\sqrt {2} \quad .
$$
Initial conditions 
$$
a_0=0 \,, \, a_1=12 \quad ; \quad b_0=1 \, , \, b_1=5 \quad .
$$

{\bf Integer-ating factor} $E(n)=lcm(1 \dots n)^3$, $\nu=3$.

$$
\delta= \frac{\log \beta -\nu}{\log \alpha + \nu} =
\frac{\log \beta -3}{\log \alpha + 3} =
\frac{\log ( 17+12\,\sqrt {2}) -3}{\log ( 17+12\,\sqrt {2}) + 3} =
0.080529431189061685186
\dots \quad.
$$

Implied irrationality measure: $1+1/\delta= 13.41782023335376578458\dots$.

{\bf Accompanying Maple packages}

This article is accompanied by three Maple packages, {\tt GenBeukersLog.txt},  {\tt GenBeukersZeta2.txt}, {\tt GenBeukersZeta3.txt}
all freely available from the {\it front} of this masterpiece

{\tt https://sites.math.rutgers.edu/\~{}zeilberg/mamarim/mamarimhtml/beukers.html} \quad ,

where one can find ample sample input and output files, that readers are welcome to extend.

{\bf Zudilin's Tweak of the Beukers $\zeta(2)$ integral to get the Catalan constant}

The {\it inspiration} for our {\it tweaks} came from Wadim Zudilin's  brilliant discovery [Zu1] that the famous Catalan constant,
that may be defined by the innocent-looking alternating series of the reciprocals of the odd perfect-squares
$$
C:=1- \frac{1}{3^2} + \frac{1}{5^2} - \frac{1}{7^2} + \dots = \sum_{n=0}^{\infty} \frac{(-1)^n}{(2n+1)^2} \quad,
$$

can be written as the double integral
$$
\frac{1}{8} \, \int_0^1 \, \int_0^1 \, \frac{x^{-\frac{1}{2}} (1-y)^{-\frac{1}{2}}}
{1-xy} \, dx \, dy \quad .
$$
This lead him to consider the sequence of Beukers-type double-integrals
$$
I(n) \, = \,
\int_0^1 \, \int_0^1 \, 
\frac{
x^{-\frac{1}{2}} 
(1-y)^{-\frac{1}{2}}
     }
{1-xy} 
\cdot \left( \frac{x(1-x)y(1-y)}{1-xy} \right )^n \, dx \, dy \quad .
$$

Using the Zeilberger algorithm, Zudilin derived a three term recurrence for $I(n)$ leading to good diophantine approximations
to the Catalan constant, alas not good enough to prove irrationality. This was elaborated and extended by Yu. V. Nesterenko [N].
See also [Zu2].

Using the multivariable Almkvist-Zeilberger algorithm we can derive the recurrence much faster. Using Koutschan's package [K], it is yet faster.

{\bf Our Tweaks}

Inspired by Zudilin's Beukers-like integral for the Catalan constant, we decided to use our efficient tools for quickly manufacturing
recurrences.

We systematically investigated the following families.

{\bf Generalizing the Alladi-Robinson-Like Integral for $\log 2$}

Alladi and Robinson [AR] gave a Beukers-style new proof of the irrationality of $\log2 $ using the
elementary fact that

$$
\log 2 \, = \, \int_0^1 \, \frac{1}{1+x} \, dx \quad ,
$$
and more generally,
$$
\frac{1}{c} \, \log (1+c) \, = \, \int_0^1 \, \frac{1}{1+cx} \, dx \quad.
$$

They used the sequence of integrals
$$
I(n):=\int_0^1 \, \frac{1}{1+cx} \left (\frac{x(1-x)}{1+cx} \right )^n \, dx \quad,
$$
and proved that for a wide range of choices of rational $c$, this leads to irrationality proofs and irrationality measures (see also [ZeZu1]).

Our generalized version is the three-parameter family of constants
$$
I_1(a,b,c):= \frac{1}{B(1+a,1+b)} \, \int_0^1 \, \frac{x^a(1-x)^b}{1+cx} \, dx
$$
that is easily  seen to  equal  ${}_2F_1(1,a+1;a+b+2;-c)$.

We use the sequence of integrals
$$
I_1(a,b,c)(n):= \,  \frac{1}{B(1+a,1+b)} \, \int_0^1 \, \frac{x^a(1-x)^b}{1+cx} \cdot \left( \frac{x(1-x)}{1+cx} \right )^n\, dx \quad .
$$

Using the (original!) Almkvist-Zeilberger algorithm [AlZ], implemented in the Maple package

{\tt https://sites.math.rutgers.edu/\~{}zeilberg/tokhniot/EKHAD.txt} \quad,

we immediately get a second-order recurrence that can be gotten by typing `{\tt OpeL(a,b,c,n,N);}' in the Maple package

{\tt https://sites.math.rutgers.edu/\~{}zeilberg/tokhniot/GenBeukersLog.txt} \quad.

This enabled us to conduct a systematic search, and we found many cases of ${}_2F_1$ evaluations that lead to
irrationality proofs, i.e. for which the $\delta$ mentioned above is {\bf positive}. Many of them
turned out to be (conjecturally) expressible in terms of algebraic numbers and/or logarithms of rational numbers,
hence proving them irrational is not that exciting, but we have quite a few not-yet-identified (and inequivalent) cases. See the output file

{\tt https://sites.math.rutgers.edu/\~{}zeilberg/tokhniot/oGenBeukersLog1.txt} \quad,

for many examples. Whenever Maple was able to (conjecturally) identify the constants explicitly, it is mentioned.
If nothing is mentioned then these are potentially explicit constants, expressible as a hypergeometric series ${}_2 F_1$,
for which this would be the {\bf first} irrationality proof, once the details are filled-in.

We also considered the four-parameter family of constants
$$
I'_1(a,b,c,d):= 
\frac
{\int_0^1 \, \frac{x^a(1-x)^b}{(1+cx)^{d+1}} \, dx}
{\int_0^1 \, \frac{x^a(1-x)^b}{(1+cx)^d} \, dx} \quad,
$$
and, using the  more general recurrence, also obtained using the Almkvist-Zeilberger algorithm
(to see it type `{\tt OpeLg(a,b,c,d,n,Sn);}' in  {\tt GenBeukersLog.txt}),
found many candidates for irrationality proofs that Maple was unable to identify. See the output file

{\tt https://sites.math.rutgers.edu/\~{}zeilberg/tokhniot/oGenBeukersLog2.txt} \quad.

{\bf Generalizing the Beukers Integral for $\zeta(2)$}

Define
$$
I_2(a_1,a_2,b_1,b_2)(n)
\, = \,
\frac{1}{B(1-a_1,1-a_2) B(1-b_1,1-b_2)} \cdot
$$
$$
\int_0^1 \, \int_0^1 \, \frac{x^{-a_1} (1-x)^{-a_2} y^{-b_1} (1-y)^{-b_2} }{1-xy} \cdot \left( \frac{x(1-x)y(1-y)}{1-xy} \right )^n \, dx \, dy \quad ,
$$
that happens to satisfy a linear-recurrence equation of second order, yielding Diophantine approximations to the constant
$I_2(a_1,a_2,b_1,b_2)(0)$, let's call it $C_2(a_1,a_2,b_1,b_2)$
$$
C_2(a_1,a_2,b_1,b_2)
\, = \,
\frac{1}{B(1-a_1,1-a_2) B(1-b_1,1-b_2)} \cdot
\int_0^1 \, \int_0^1 \, \frac{x^{-a_1} (1-x)^{-a_2} y^{-b_1} (1-y)^{-b2}}{1-xy} \, dx \, dy \quad .
$$

It is readily seen that 
$$
C_2(a_1,a_2,b_1,b_2) \, = \,{}_3 F_2
\left (
{{1 \,, \, 1 - a_1 \, , \, -b_1 + 1}
\atop
{ 2 - a_1 - a_2 \, , \, 2 - b_1 - b_2}} \, ;1 \, \right) \quad .
$$

Most choices of  random $a_1,a_2,b_1,b_2$ yield disappointing, negative $\delta$'s, just like $C_2(\frac{1}{2},0,0,\frac{1}{2})$ (alias  $8$ times the
Catalan constant), but a systematic search yielded several hundred candidates
that produce positive $\delta$'s and hence would produce irrationality proofs.
Alas, many of them were conjecturally equivalent to each other
via  a fractional-linear transformation with integer coefficients, $C \rightarrow \frac{a+bC}{c+dC}$, with $a,b,c,d$ integers, hence
the facts that they are irrational are equivalent. Nevertheless we found quite a few that are (conjecturally) {\it not equivalent} to each other.
Modulo filling-in some details, they lead to irrationality proofs. Amongst them some were (conjecturally) identified by
Maple to be either algebraic, or logarithms of rational numbers, 
for which irrationality proofs exist for thousands of years (in case of $\sqrt{2}$ and $\sqrt{3}$ etc.),
or a few hundred years (in case of $\log 2$, etc.).

But some of them Maple was {\bf unable} to identify, so {\it potentially} our (sketches) of proofs would be the {\bf first} irrationality proofs.

{\bf Beukers $\zeta(2)$  Tweaks That produced Irrationality Proofs with Identified Constants}

{\bf Denominator 2}

We first searched for $C_2(a_1,a_2,b_1,b_2)$ where the parameters $a_1,a_2,b_1,b_2$ have  denominator $2$, 
there were quite a few of them, but they were
all conjecturally equivalent to each other. Here is one of them:

$\bullet$  $C_2(0,0,\frac{1}{2},0)={}_3F_2(1,1,1/2;2,3/2;1)$, alias $2 \log 2$.

{\bf Denominator 3}

There were also quite a few where the parameters $a_1,a_2,b_1,b_2$ have  denominator $3$,
but again they were all equivalent to each other, featuring
$\pi \sqrt{3}$. Here is one of them.

$\bullet$  $C_2(0,0,\frac{1}{3},-\frac{2}{3}) = {}_3F_2(1,1,2/3;2,7/3;1)$, alias (conjecturally) $-6+4\pi\sqrt{3}/3$.

{\bf Denominator 4}

There were also quite a few 
where the parameters $a_1,a_2,b_1,b_2$ have  denominator $4$, 
but again they were all equivalent to each other, featuring
$\sqrt{2}$, yielding a new proof of the irrationality of $\sqrt{2}$ (for what it is worth). Here is one of them.

$\bullet$  $C_2(-\frac{3}{4}, -\frac{3}{4}, -\frac{1}{4}, -\frac{3}{4})={}_3F_2(1,7/4,5/4; 7/2,3;1)$, alias (conjecturally) $-240\,+\,\frac{512}{3}\,\sqrt{2}$.

{\bf Denominator 5}

There were also quite a few 
where the parameters $a_1,a_2,b_1,b_2$ have  denominator $5$, but again they were all equivalent to each other, featuring
$\sqrt{5}$, yielding a new proof of the irrationality of $\sqrt{5}$ (for what it is worth). Here is one of them.

$\bullet$  $C_2(-\frac{4}{5}, -\frac{4}{5},-\frac{2}{5},-\frac{3}{5} )={}_3F_2(1,9/5,7/5; 18/5,3;1)$, 
alias (conjecturally) $-\frac{845}{2}\, + \,\frac{2275}{12}\,\sqrt{5}$

{\bf Denominator 6 with identified constants}

We found two equivalence classes 
where the parameters $a_1,a_2,b_1,b_2$ have  denominator $6$, for which the constants were identified.
Here are one from each class.

$\bullet$  $C_2( -5/6, -5/6, -1/2, -1/2)={}_3F_2(1,11/6,3/2; 11/3,3;1)$, alias (conjecturally) $-{\frac{1344}{5}}+{\frac {16384\,\sqrt {3}}{105}}$

$\bullet$  $C_2( -5/6, -5/6, -1/3, -2/3)={}_3F_2(1,11/6,4/3; 11/3,3;1)$, alias (conjecturally) ${\frac {972\,{2}^{2/3}}{5}}-{\frac{1536}{5}}$

{\bf denominator 7 with identified constants}

We found two cases 
where the parameters $a_1,a_2,b_1,b_2$ have  denominator $7$, for which the constants were identified.

$\bullet$  $C_2( -6/7, -6/7, -4/7, -3/7)={}_3F_2(1,13/7,11/7;26/7,3;1)$, alias (conjecturally) the positive root of $13824\,{x}^{3}-2757888\,{x}^{2}-10737789048\,x+16108505539=0$ .

$\bullet$  $C_2( -6/7, -1/7, 4/7, 2/7)={}_3F_2( 1,13/7,3/7;3,8/7;1)$, alias (conjecturally) the positive root of $2299968\,{x}^{3}+7074144\,{x}^{2}-11234916\,x-12663217=0$

{\bf Beukers $\zeta(2)$  Tweaks That produced Irrationality Proofs with Not-Yet-Identified Constants (and Hence Candidates for First Irrationality Proofs)}

For the following constants, Maple was unable to identify, and we have {\it potentially} the first irrationality proofs of these
constants.

{\bf Denominator 6 with not yet identified constants}

We found two cases (up to equivalence):

$\bullet$ $C_2( 0, -1/2, 1/6, -1/2)={}_3F_2(1,1,5/6;5/2,7/3;1)$

While Maple was unable to identify this constant, Mathematica came up with $-24 \,-\, \frac{81 \sqrt{\pi} \Gamma(7/3)}{\Gamma(-1/6)}$.

$\bullet$ $C_2( -2/3, -1/2, 1/2, -1/2)={}_3F_2(1,5/3,1/2;19/6,2;1)$

While Maple was unable to identify this constant, Mathematica came up with $\frac{13}{2}\, -\,\frac{6 \Gamma(19/6)}{ \sqrt{\pi}\Gamma(8/3)}$.

{\bf Denominator 7 with not yet identified constants}

We found six cases (up to equivalence):

$\bullet$ $C_2( -6/7, -6/7, -4/7, -5/7)={}_3F_2( 1,13/7,11/7;26/7,23/7;1)$

$\bullet$ $C_2( -6/7, -5/7, -3/7, -5/7)={}_3F_2( 1,13/7,10/7;25/7,22/7;1)$

$\bullet$ $C_2( -6/7, -5/7, -2/7, -1/7)={}_3F_2( 1,13/7,9/7;25/7,17/7;1)$

$\bullet$ $C_2( -6/7, -4/7, -1/7, -1/7)={}_3F_2( 1,13/7,8/7;24/7,16/7;1)$

$\bullet$ $C_2( -6/7, -3/7, -5/7, -3/7)={}_3F_2( 1,13/7,12/7;23/7,22/7;1)$

$\bullet$ $C_2( -5/7, -3/7, -4/7, -2/7)={}_3F_2( 1,12/7,11/7;22/7,20/7;1)$

For each of them, to get the corresponding theorem and proof, use procedure {\tt TheoremZ2} in the Maple pacgage
{\tt GenBeukersZeta2.txt}.

To get a statement and full proof (modulo a divisibility lemma) type , in {\tt GenBeukersZeta2.txt}

{\tt TheoremZ2(a1,a2,b1,b2,K,0):}

with {\tt K} at least $2000$. For example, for the last constant in the above list ${}_3F_2( 1,12/7,11/7;22/7,20/7;1)$, type

{\tt TheoremZ2( -5/7, -3/7, -4/7, -2/7 ,3000,0):}

For more details (the recurrences, the estimated irrationality measures, the initial conditions)
see the output file

{\tt https://sites.math.rutgers.edu/\~{}zeilberg/tokhniot/oGenBeukersZeta2g.txt} \quad .

{\bf Generalizing the Beukers Integral for $\zeta(3)$}

The natural extension would be the {\bf six}-parameter family (but now we make the exponents positive)
$$
\frac{1}{B(1+a_1,1+a_2) B(1+b_1,1+b_2) B(1+c_1,1+c_2)} \cdot
$$
$$
\int_0^1 \, \int_0^1 \, \int_0^1 \, 
\frac{x^{a_1} (1-x)^{a_2} y^{b_1} (1-y)^{b_2} z^{c_1} (1-z)^{c_2}}
{1-z+xyz} \cdot \left( \frac{x(1-x)y(1-y)z(1-z)}{1-z+xyz} \right )^n \, dx \, dy \, dz \quad .
$$

However, for {\bf arbitrary} $a_1,a_2,b_1,b_2,c_1,c_2$ the recurrence is {\bf third order}.
(Wadim Zudilin pointed out that this may be related to the work of Rhin and Viola in [RV]).

Also, empirically, we did not find many promising cases. Instead, let's define
$$
J_3(a_1,a_2,b_1,b_2,c_1,c_2;e)(n)
$$
$$
\int_0^1 \, \int_0^1 \, \int_0^1 \, 
\frac{x^{a_1} (1-x)^{a_2} y^{b_1} (1-y)^{b_2} z^{c_1} (1-z)^{c_2}}{(1-z+xyz)^e} 
\cdot \left( \frac{x(1-x)y(1-y)z(1-z)}{1-z+xyz} \right )^n \, dx \, dy \, dz\quad .
$$
and
$$
I_3(a_1,a_2,b_1,b_2,c_1,c_2;e)(n) :=
\frac{J_3(a_1,a_2,b_1,b_2,c_1,c_2;e+1)(n)}{J_3(a_1,a_2,b_1,b_2,c_1,c_2;e)(0)}
$$

The family of constants that we hope to prove irrationality is the five-parameter:
$$
I_3(a_1,a_2,b_1,b_2,c_1,c_2;e)(0) \quad.
$$
$$
=\frac
{
\int_0^1 \, \int_0^1 \, \int_0^1 \, \frac{x^{a_1} (1-x)^{a_2} y^{b_1} (1-y)^{b_2} z^{c_1} (1-z)^{c_2}}{(1-z+xyz)^{e+1}} \,dx\,dy \, dz
}
{
\int_0^1 \, \int_0^1 \, \int_0^1 \, \frac{x^{a_1} (1-x)^{a_2} y^{b_1} (1-y)^{b_2} z^{c_1} (1-z)^{c_2}}{(1-z+xyz)^e} \,dx\,dy \, dz
} \quad.
$$
Of course, for  this more general, $7$-parameter, family, there is no second-order recurrence, but rather
a third-order one. But to our delight, we found a five-parameter family, let's call it
$$
K(a,b,c,d,e)(n):= I_3(b, c, e, a, a, c, d)(n) \quad.
$$
Spelled-out, our five-parameter family of constants is
$$
K(a,b,c,d,e)(0)=
$$
$$
\frac
{
\int_0^1 \, \int_0^1 \, \int_0^1 \, \frac{x^{b} (1-x)^{c} y^{e} (1-y)^{a} z^{a} (1-z)^{c}}{(1-z+xyz)^{d+1}} \,dx\,dy \, dz
}
{
\int_0^1 \, \int_0^1 \, \int_0^1 \, \frac{x^{b} (1-x)^{c} y^{e} (1-y)^{a} z^{a} (1-z)^{c}}{(1-z+xyz)^{d}} \,dx\,dy \, dz
} \quad.
$$

Now we found (see the section on finding recurrences below) a general second-order recurrence, that  is too complicated
to display here in full generality, but can be seen by typing

{\tt OPEZ3(a,b,c,d,e,n,Sn);}

In the Maple package {\tt GenBeukersZeta3.txt}. This enabled us, for each specific, numeric specialization of
the parameters $a,b,c,d,e$ to quickly find the relevant recurrence, and  systematically search for
those that give positive $\delta$. Once again, many of them turned out to be (conjecturally) equivalent to
each other.

{\bf Denominator 2}:

We only found one class, up to equivalence, all related to $\log 2$. One of them is
$$
K(0,0,0,1/2,1/2)=I_3(0, 0, 1/2, 0, 0, 0, 1/2) \quad,
$$
that is not that exciting since it is (conjecturally) equal to $-{\frac {2-4\,\log  \left( 2 \right) }{3-4\,\log  \left( 2 \right) }}$.

For details, type {\tt TheoremZ3(0,0,0,1/2,1/2,3000,0);} in {\tt GenBeukersZeta3.txt} \quad .

\vfill\eject

{\bf Denominator 3}:

We found three inequivalent classes, none of them Maple was able to identify.

$$
K( 0, 0, 0, 1/3,2/3)=I_3( 0, 0, 2/3, 0, 0, 0, 1/3) \quad,
$$
for details, type {\tt TheoremZ3(0,0,0,1/3,2/3,3000,0);} in {\tt GenBeukersZeta3.txt}.

$$
K(0,0,0,2/3,1/3 )=I_3( 0, 0, 1/3, 0, 0, 0, 2/3) \quad,
$$
for details, type {\tt TheoremZ3(0,0,0,2/3,1/3,3000,0);} in {\tt GenBeukersZeta3.txt}.

$$
K(0, 1/3, 2/3, 1/3, 2/3 )=I_3(  0, 0, 1/3, 0, 0, 0, 2/3) \quad,
$$
for details, type {\tt TheoremZ3(0,1/3,2/3,1/3,2/3,3000,0);} in {\tt GenBeukersZeta3.txt},

These three constants are {\bf candidates} for  {\bf `first-ever-irrationality proof'}.

{\bf Denominator 4}: We only found one family, all expressible in terms of $\log 2$. Here is one of them.

For example
$$
K( 0, 1/2, 0, 1/4, 3/4)= I_3(1/2, 0, 3/4, 0, 0, 0, 1/4) \quad,
$$
that, conjecturally  equals $-{\frac {-30+45\,\log  \left( 2 \right) }{-11+15\,\log  \left( 2 \right) }}$.

For details, type {\tt TheoremZ3(0,1/2,0,1/4,3/4,3000,0);} in {\tt GenBeukersZeta3.txt}.

{\bf Denominator 5}: We only found one family, up to equivalence, but Maple was unable to identify the constant. So it is
potentially the first irrationality proof of that constant
$$
K(0, 1/5, 0, 3/5, 2/5)=I_3(1/5, 0, 2/5, 0, 0, 0, 3/5) \quad .
$$
For details, type {\tt TheoremZ3(0,1/5,0,3/5,2/5,3000,0);} in {\tt GenBeukersZeta3.txt}.

{\bf Denominator 6}: We  found three families, up to equivalence,  none of which Maple was able to identify. Once again,
these are candidates for {\bf first-ever irrationality proofs} for these constants.

$$
K(0, 1/2, 1/2, 1/3, 1/6)=I_3(1/2, 1/2, 1/6, 0, 0, 1/2, 1/3) \quad .
$$
For details, type {\tt TheoremZ3(0,1/2,1/2,1/3,1/6,3000,0);} in {\tt GenBeukersZeta3.txt}.

$$
K(0, 1/2, 1/2, 1/6, 1/3)=I_3(1/2, 1/2, 1/3, 0, 0, 1/2, 1/6) \quad.
$$
For details, type {\tt TheoremZ3(0,1/2,1/2,1/6,1/3,3000,0);} in {\tt GenBeukersZeta3.txt}.

$$
K(1/3, 0, 2/3, 1/2, 5/6)=I_3(0, 2/3, 5/6, 1/3, 1/3, 2/3, 1/2) \quad .
$$
For details, type {\tt TheoremZ3(1/3,0,2/3,1/2,5/6,3000,0);} in {\tt GenBeukersZeta3.txt}.

{\bf Denominator 7}: We  found five families, up to equivalence,  none of which Maple was able to identify. Once again,
these are candidates for first-ever irrationality proofs for these constants.

$$
K(1/7, 0, 2/7, 3/7, 4/7)=I_3(0, 2/7, 4/7, 1/7, 1/7, 2/7, 3/7) \quad .
$$
For details, type {\tt TheoremZ3(1/7,0,2/7,3/7,4/7,3000,0);} in {\tt GenBeukersZeta3.txt}.

$$
K( 1/7, 0, 2/7, 5/7, 3/7)=I_3(0, 2/7, 3/7, 1/7, 1/7, 2/7, 5/7) \quad .
$$
For details, type {\tt TheoremZ3(1/7,0,2/7,5/7,3/7,3000,0);} in {\tt GenBeukersZeta3.txt}.

$$
K(1/7, 0, 3/7, 4/7, 5/7)=I_3(   0, 3/7, 5/7, 1/7, 1/7, 3/7, 4/7) \quad .
$$
For details, type {\tt TheoremZ3(1/7,0,3/7,4/7,5/7,3000,0);} in {\tt GenBeukersZeta3.txt}.

$$
K(1/7, 0, 4/7, 2/7, 5/7)=I_3(0, 4/7, 5/7, 1/7, 1/7, 4/7, 2/7) \quad .
$$
For details, type {\tt TheoremZ3(1/7,0,4/7,2/7,5/7,3000,0);} in {\tt GenBeukersZeta3.txt}.

$$
K( 2/7, 0, 3/7, 4/7, 5/7)=I_3(0, 3/7, 5/7, 2/7, 2/7, 3/7, 4/7) \quad .
$$
For details, type {\tt TheoremZ3(2/7,0,3/7,4/7,5/7,3000,0);} in {\tt GenBeukersZeta3.txt}.

If you don't have Maple, you can look at the output file

{\tt https://sites.math.rutgers.edu/\~{}zeilberg/tokhniot/oGenBeukersZeta3All.txt} \quad ,

that gives detailed sketches of irrationality proofs of all the above constants, some with conjectured
integer-ating factors.

{\bf Guessing an INTEGER-ating factor}

In the original Beukers cases the integer-ating factor was easy to conjecture, and even to prove.
For $\zeta(2)$ it was $lcm(1 \dots n)^2$, and for $\zeta(3)$ it  was  $lcm(1 \dots n)^3$.
For the Alladi-Robinson case of $\log 2$ it was even simpler, $lcm(1 \dots n)$.

But in other cases it is much more complicated. A natural `atomic' object is, given
a modulo M, a subset C of $\{0,..., M-1\}$, rational numbers $e_1$, $e_2$ between $0$ and $1$, rational numbers $e_3,e_4$,
the following quantity, for positive integers $n$
$$
Pp(e_1,e_2,e_3,e_4, C, M; n):= \prod_p p \quad,
$$
where $p$ ranges over all primes such that (let $\{a\}$ be the fractional part of $a$, i.e. $a- \lfloor a \rfloor$)

$\bullet$ $ e_1 < \{ n/p \} <e_2$

$\bullet$ $ e_3 < p/n <e_4$

$\bullet$ $p \,\, mod \,\, M \in C$

Using the prime number theorem, it follows (see e.g. [Zu2]) that
$$
\lim_{n \rightarrow \infty} \frac{\log Pp(e_1,e_2,e_3,e_4, C, M; n)}{n} \quad,
$$
can be evaluated exactly, in terms of the function $\Psi(x)=\frac{\Gamma'(x)}{\Gamma(x)}$ (see procedure {\tt PpGlimit} in the Maple packages)
thereby giving an exact value for the quantity $\delta$ whose positivity implies irrationality.

Of course, one still needs to rigorously prove that the conjectured integer-ating factor is indeed correct.

{\bf Looking under the hood: On Recurrence Equations}

For `secrets from the kitchen' on how we found the second-order, four-parameter recurrence operator \hfill\break
{\tt OPEZ2(a1,a2,b1,b2,n,N)} 
in the Maple package {\tt GenBeukersZeta2.txt},
that was the engine driving the $\zeta(2)$ tweaks, and more impressively,
the five-parameter second-order recurrence operator {\tt OPEZ3(a,b,c,d,e,n,N)}
in the Maple package {\tt GenBeukersZeta3.txt}, that was the engine driving the $\zeta(3)$ tweaks,
the reader is referred to the stand-alone appendix available from the following url:

{\tt https://sites.math.rutgers.edu/\~{}zeilberg/mamarim/mamarimPDF/beukersAppendix.pdf} \quad .

{\bf Other Variations on Ap\'ery's theme}

Other attempts to use Ap\'ery's brilliant insight are [Ze2][Ze3][ZeZu1]. Recently Marc Chamberland and Armin Straub [CS]
explored other fascinating aspects of the Ap\'ery numbers, not related to irrationality.

{\bf Conclusion and Future Work}

We believe that symbolic computational methods  have great potential in irrationality proofs,
in particular, and number theory in general.
In this article we confined attention to approximating sequences that arise from {\it second-order recurrences}.
The problem with higher order recurrences is that one gets linear combinations with rational coefficients of {\it several} constants,
but if you can get two different such sequences coming from third-order recurrences,
both featuring the same two constants, then the present method may be applicable.
More generally if you have a $k$-th order recurrences, you need $k-1$ different integrals. 

The general methodology of this article can be called {\it Combinatorial Number Theory}, but not in the usual sense, but
rather as an analog of {\it Combinatorial Chemistry}, where one tries out many potential chemical compounds, most of them
useless, but since computers are so fast, we can afford to generate lots of cases and pick the wheat from the chaff.

{\bf  Encore: Hypergeometric challenges}

As a {\it tangent}, we (or rather Maple) discovered many exact ${}_3F_2(1)$ evaluations.
Recall that the Zeilberger algorithm can prove hypergemoetric identities only if there is at least one free parameter.
For a {\it specific} ${}_3 F_2(a_1 \, a_2 \, a_3 \, ; b_1 \, b_2 ; 1)$, with {\it numeric} parameters, it is useless.
Of course, it is sometimes possible to introduce such a parameter in order to conjecture a {\it general} identity, valid for `infinitely' many $n$,
and then specialize $n$ to a specific value, but this remains an art rather than a science. The output file

{\tt https://sites.math.rutgers.edu/\~{}zeilberg/tokhniot/oGenBeukersZeta2f.txt} \quad 

contains many such conjectured evaluations, (very possibly many of them are equivalent via a hypergeometric
transformation rule) and we challenge Wadim Zudilin, the {\it birthday boy}, or anyone else, to prove them.

{\bf References}

[AR] Krishna Alladi and Michael L. Robinson,
{\it Legendre polynomials and irrationality},
J. Reine Angew. Math. {\bf 318} (1980), 137-155.

[AlZ] Gert Almkvist and Doron Zeilberger, {\it The method of differentiating under the
integral sign}, J. Symbolic Computation {\bf 10}, 571-591 (1990). \hfill\break
{\tt https://sites.math.rutgers.edu/\~{}zeilberg/mamarim/mamarimhtml/duis.html} \quad .

[ApaZ] Moa Apagodu  and Doron Zeilberger,
{\it Multi-variable Zeilberger and Almkvist-Zeilberger algorithms and the
sharpening of Wilf-Zeilberger Theory },
Adv. Appl. Math. {\bf 37} (2006)(Special Regev issue), 139-152. \hfill\break
{\tt https://sites.math.rutgers.edu/\~{}zeilberg/mamarim/mamarimhtml/multiZ.html} \quad .

[Ape] Roger Ap\'ery,
{\it ``Interpolation de fractions continues et irrationalit\'e
de certaine constantes''}
Bulletin de la section des sciences du C.T.H.S. \#3
p. 37-53, 1981.

[B] Frits Beukers, {\it A note on the irrationality of
$\zeta(2)$ and $\zeta(3)$}, 
Bull. London Math. Soc. {\bf 11} (1979), 268-272.

[CS]  Marc Chamberland and Armin Straub, {\it Ap\'ery limits: Experiments and Proofs}, arxiv:2001.034400v1, 6 Nov 2020. \hfill\break
{\tt https://arxiv.org/abs/2011.03400} \quad .

[H] Professor David Hilbert, {\it Mathematical Problems} [Lecture delivered before the International Congress of Mathematicians at
Paris in 1900], translated by Dr. Mary Winston Newson, Bulletin of the American Mathematica Society {\bf 8} (1902), 437-479. \hfill\break
{\tt https://www.ams.org/journals/bull/2000-37-04/S0273-0979-00-00881-8/S0273-0979-00-00881-8.pdf} \quad .

[K] Christoph Koutschan, {\it   Advanced applications of the holonomic systems approach}, 
PhD thesis, Research Institute for Symbolic Computation (RISC), Johannes Kepler University, Linz, Austria, 2009.\hfill\break
{\tt http://www.koutschan.de/publ/Koutschan09/thesisKoutschan.pdf}, \hfill\break
{\tt http://www.risc.jku.at/research/combinat/software/HolonomicFunctions/} \quad.

[N] Yu. V. Nesterenko, {\it On Catalan's constant}, Proceedings of the Steklov Institute of Mathematics {\bf 292} (2016), 153-170.

[vdP] Alf van der Poorten, {\it A proof that Euler missed...
Ap\'ery's proof of the irrationality of $\zeta(3)$},
Math. Intelligencer {\bf 1} (1979), 195-203.

[RV] Georges Rhin and Carlo Viola, {\it The group structure of $\zeta(3)$}, Acta Arithmetica, {\bf 97}(2001), 269-293.

[Ze1] Doron Zeilberger, {\it A Holonomic systems approach to special functions
identities}, J. of Computational and Applied Math. {\bf 32}, 321-368 (1990). \hfill\break
{\tt  https://sites.math.rutgers.edu/\~{}zeilberg/mamarim/mamarimhtml/holonomic.html }\quad .

[Ze2] Doron Zeilberger, {\it Computerized deconstruction}, Adv. Applied Math. {\bf 30} (2003), 633-654. \hfill\break
{\tt https://sites.math.rutgers.edu/\~{}zeilberg/mamarim/mamarimhtml/derrida.html} \quad  .

[Ze3] Doron Zeilberger, {\it Searching for Ap\'ery-style miracles
[using, inter-alia, the amazing Almkvist-Zeilberger algorithm]}, Personal Journal of Shalosh B. Ekhad and Doron Zeilberger, \hfill\break
{\tt https://sites.math.rutgers.edu/\~{}zeilberg/mamarim/mamarimhtml/apery.html} \quad .

[ZeZu1] Doron Zeilberger, and Wadim Zudilin,
{\it Automatic discovery of irrationality proofs and irrationality measures},
International Journal of Number Theory , published on-line before print, volume and page tbd.
Also to appear in a book dedicated to Bruce Berndt. \hfill\break
{\tt https://sites.math.rutgers.edu/\~{}zeilberg/mamarim/mamarimhtml/gat.html} \quad .

[ZeZu2] Doron Zeilberger, and Wadim Zudilin,
{\it The irrationality measure of Pi is at most 7.103205334137...},
 Moscow J. of Combinatorics and Number Theory 9 (2020), 407-419. \hfill\break
{\tt https://sites.math.rutgers.edu/\~{}zeilberg/mamarim/mamarimhtml/pimeas.html} \quad .

[Zu1] Wadim Zudilin, {\it Ap\'ery-like difference equations for Catalan's constant}  \hfill\break
{\tt https://arxiv.org/abs/math/0201024} \quad .

[Zu2] Wadim Zudilin, {\it  Arithmetic of linear forms involving odd zeta values},
      J. Th\'eorie Nombres Bordeaux {\bf 16} (2004), 251-291. \hfill\break
{\tt https://arxiv.org/abs/math/0206176} \quad .

\bigskip
\hrule
\bigskip
Robert Dougherty-Bliss, Department of Mathematics, Rutgers University (New Brunswick), Hill Center-Busch Campus, 110 Frelinghuysen
Rd., Piscataway, NJ 08854-8019, USA. \hfill\break
Email: {\tt Robert.w.Bliss at gmail  dot com}   \quad .
\bigskip
Christoph Koutschan, Johann Radon Institute of Computational and Applied Mathematics (RICAM), Austrian Academy of Sciences,
Altenberger Strasse 69, A-4040 Linz, Austria \hfill\break
Email: {\tt  christoph.koutschan at ricam dot oeaw dot ac dot at}   \quad .
\bigskip
Doron Zeilberger, Department of Mathematics, Rutgers University (New Brunswick), Hill Center-Busch Campus, 110 Frelinghuysen
Rd., Piscataway, NJ 08854-8019, USA. \hfill\break
Email: {\tt DoronZeil at gmail  dot com}   \quad .

\end